\documentclass[11pt,a4paper]{article}
\usepackage[utf8]{inputenc}
\usepackage[english]{babel}
\usepackage{amssymb, amsmath, enumerate}
\usepackage{graphicx}
\usepackage{epstopdf}
\usepackage{subcaption}
\usepackage{authblk}
\usepackage{textcomp}

\newcommand{\qed}{\nobreak \ifvmode \relax \else
	\ifdim\lastskip<1.5em \hskip-\lastskip
	\hskip1.5em plus0em minus0.5em \fi \nobreak
	\vrule height0.75em width0.5em depth0.25em\fi}
\begin{document}
	\title{ Preventive Equipment Repair Planning Model}
	\author[1]{O. A. Malafeyev\thanks{o.malafeev@spbu.ru}}
	\author[1]{N.D.Redinskikh\thanks{redinskich@yandex.ru}}
	\author[1]{S. A. Nemnyugin\thanks{s.nemnyugin@spbu.ru}}
	\author[1]{I. D. Kolesin\thanks{kolesin$_{-}$id@mail.ru}}
	\author[2,3]{I. V. Zaitseva\thanks{ziki@mail.ru}}
	\affil[1]{St.Petersburg State University, 7/9 Universitetskayanab., St. Petersburg, 199034 Russia}
	\affil[2]{Stavropol State Agrarian University, Zootekhnicheskiy lane 12, Stavropol, 355017, Russia.}
	\affil[3]{Stavropol branch of the Moscow Pedagogical State University, Dovatortsev str. 66 g, Stavropol, 355042, Russia.}
	\date{}
	\maketitle
	\begin{abstract}
	Any large functioning system consists of equipment that needs to be repaired during its lifetime. The methods of mathematical programming are used to formalize the optimization problem of preventive equipment repair planning in this paper. The method is used to find the annual optimal plan of preventive equipment repair. The goal of such plan is the uniform distribution of repair works by the months. It is supposed, that the deviation from the standard plan of repair works can not be more than one month. An algorithm to find the optimal plan of repair works is provided.  A numerical example is given. 
	\end{abstract}
	
	

	\section{Introduction}
	
	Different problems of preventive equipment repair using various mathematical methods are considered in many papers. For example, in \cite{HMWagner} the methods of mathematical programming are applied to study the problems of preventive maintenance scheduling. In \cite{HMWagner} authors illustrate how linear programming can be effectively used to solve problems of such class and  propose different optimality criteria. Several articles discuss the problem of preventive maintenance efficacy testing. The hypothesis that repairs have no effect on the sequence of times to failure or on the costs of the failures is tested in \cite{RDBaker}. Other mathematical models that can be used for formalization and analysis of the problem of preventive repair planning are illustrated in [3--51].
	
	Equipment which fails during its lifetime and needs to be repaired is considered in this paper. For such cases, it is usually necessary to ensure a uniform plan of repair works. This plan should contain information about the equipment repair periods for a year. Such plan can be considered as a normative annual plan of equipment repair. It is formed in accordance with standard monthly volumes of the repair works. But usually, in real life, the volumes of repair works do not correspond to the actual plan of equipment repairs.
	
	Therefore there is a need of regulatory adjustments to the annual plan of equipment repairs through changes in the time periods of repair works. 
	However, such changes in the time periods of repair works should not be great.
	It is supposed that repair works can be moved for no more than one month (forward or backward) within a year. Under this constrain, it is possible not to obtain a uniform distribution of monthly volume of repair works. 
	Then the difference between the 
	normative annual plan of equipment repairs and a new distribution of monthly volume of repair works can be analyzed. If this difference is more than one month, then the correction of the normative annual plan can be carried out. 
	And in this case the revised plan is taken as the new normative annual plan.

\section{The formalization of the optimization problem of planning preventive equipment repair}

Let us suppose $A_1,\dots,A_{12}$ are the monthly values of the repair works expressed in repair units(hours). Each $j$-th monthly repair volume consists of the elements 
$\left(
\begin{array}{c}
a_{1j}\\
a_{2j}\\ 
\dots\\
a_{kj}
\end{array}
\right)$, 
$j=\overline{1,12}$, where $A_{ij}=\sum\limits_{i=1}^{k}a_{ij}, j=\overline{1,12}$.
Thus we obtain the matrix A=$\left(
\begin{array}{c c c}
a_{1,1}&\dots&a_{1,12}\\
a_{2,1}&\dots&a_{2,12}\\
&\dots&\\
a_{k,1}&\dots&a_{k,12}\\
\end{array}
\right)$, where $a_{ij}\geq0$, $A$ is a standard matrix of normative annual equipment repair  plan. $A_j=\sum\limits_{i=1}^{k} a_{ij}$  is the total monthly work hours in a given month $j$. Indices $ij$ give the information about equipment, which needs to be fixed this month: a) if $a_{ij}>0$, then $ij$-th equipment needs to be fixed. The repair will last $a_{ij}$ hours;
b)if $a_{ij}=0$, then $ij$-th equipment does not need to be repaired. The average monthly repair works last $A_{cp}=\frac{\sum\limits_{i=1}^{12}A_j}{12}=\frac{\sum\limits_{i=1}^k\sum\limits_{j=1}^12 a_{ij}}{12}$. The transfer of repair works can be carried out only one month backward or forward. This means, that the element $a_{ij}$ can be transferred either to the column $(j+1)$, or to the $(j-1)$-th column. It is necessary to find transfers, that have a minimal deviation the monthly volumes of repair works from the average volume of repair works. 

If $A_j^{'}$ is the new monthly volume of repair works in hours, then it is necessary to minimize the expression $\sum\limits_{j=1}^{12}|A_{j}^{'}-A_{avg}|$.
Let us denote by $x_k>0$ the volume of repair works transferred from $k$-th month to $(k+1)$-th. Let us denote by $x_k<0$ the volume of repair works transferred from $(k+1)$-th month to $k$-th. If $x_k=0$ there are no transfers. Here $x_k$ is integer number.
Let $X=(x_{ij})_{k\times 11},i=\overline{1,k},j=\overline{1,11}$ be a matrix, where
\[
\left\{
\begin{array}{l}
0 ,\small{\textnormal{if nothing is transferred}}\\
I, \small{\textnormal{if some volume of repair works is transferred from the \textit{j}-th month to \textit{(j+1)}-th}}\\
-I, \small{\textnormal{if some volume of repair works is transferred from  the \textit{(j+1)}-th month to \textit{j}-th}}\\
\end{array}
\right.
\]

Then the problem is subdivided in two problems.

\textbf{Problem 1.} It's necessary to find the monthly volumes of transferable repair works in hours, that is $x_k, k=\overline{1,11}$.

\textbf{Problem 2.} The problem of integer choice, that is the calculation of the matrix $X$.

\textit{The solution for the problem 1.}

\textit{Method I.} 
The problem can be solved by methods of quadratic programming. We use the quadratic function as the measure for deviation of the repair works volume from the average volume of repair works \begin{equation}\label{I}
V=\sum\limits_{i=1}\left(A_i-A_{avg}-x_i+x_{i-1}\right)^2+\left(A_n-A_{avg}+x_{n-1}\right)^2,
\end{equation} that should be minimized subject to \begin{equation}\label{II}
-A_{i+1}\leq x_i \leq A_i, i\div n-1
\end{equation} and additional condition of integrality $x_i$.
Consider the following transformations: $\hat{A}_i=A_i-A_{avg}$. Then the function (\ref{I}) takes the form $V=2\sum\limits_{i=1}^{n-1}\left(\hat{A}_{i+1}-\hat{A}_i\right)x_i+\sum\limits_{i=1}^{n-1}x_i^2-2\sum\limits_{i=2}^{n-1}x_{i-1}x_i+\sum\limits_{i=1}^n \hat{A}_i^2$, $\hat{A}_i=A_i-A_{avg}$. Let us suppose $z=-(V-\sum\limits_{i=1}^n \hat{A}_i^2)$. Then
 \begin{equation}\label{III}z=-\sum\limits_{i=1}^{n-1}\left(\hat{A}_{i+1}-\hat{A}_i\right)x_i+2\sum\limits_{i=2}^{n-1}x_{i-1}x_i-2\sum\limits_{i=1} x_i^2.\end{equation}.
 Let us introduce the variables $x_i^{'}$ and $x_i^{''}$ to transform the inequality (\ref{II}) into equations. Then restrictions will take the form of
\[
\left\{
\begin{array}{l}
x_i+x_i^{'}=A_i\\      
-x_i+x_i^{''}=A_{i+1}\\
x_i^{'}\geq 0\\
x_i^{''}\leq 0
\end{array}
\right.
\]

Let us introduce a new variable $x_0\geq 0$ and express each variable $x_i$ as a difference: $x_i=\overline{x_i}-x_0$, where $\overline{x_i}\geq 0$; $x_0\geq0$. Then the target function $V$ would be $V=\sum\limits_{i=1}^{n-1}\left(\hat{A}_i-(\overline{x_i}-x_0)+(\overline{x_i}-x_0)\right)^2+\left(A_n+(\overline{x_{n-1}}-x_0)-A_{avg}\right)^2$.
Let us rewrite the equation (\ref{III}) in the form 
$z=-2\sum\limits_{i=1}^{n-1}\left(\hat{A}_{i+1}-\hat{A}_i\right)\overline{x_i}-x_0\left(2\hat{A}_i-2\hat{A}_n\right)+\sum\overline{x_i}^2-2x_0^2$. Since $\min z = \max(-z)$. The problem can be represented in the form:
$$\max z = -2\sum\limits_{i=1}^{n-1}\left(\hat{A}_{i+1}-A_c\right)\overline{x_i}-x_0\left(2\hat{A}_1-2\hat{A}_n\right)+2\overline{x_1}x_0+2\overline{x_{n-1}}x_0-$$
 \begin{equation}\label{IV}
-2\sum\limits_{i=1}^{n-1}\overline{x_i}^2-2x_0^2
\end{equation} under constraints:
\begin{equation}\label{V}
(\overline{x_i}-x_0)+x_i^{'}=A_i,-(\overline{x_i}-x_0)+x_{i}^{''}=A_{i+1}, \overline{x_i}\geq 0, x_0\geq 0, x_i^{'}\geq 0, x_i^{''}\geq 0
\end{equation}  and $\overline{x_1},x_0$ are integers. The problem (\ref{IV}),(\ref{V}) is the standard problem of mathematical programming in the form:  $\max z=cx+x'Dx, A_x=B, x\geq 0$. 

\textbf{Method 2.}  The algorithm is as follows: first, the value of $x^6$ is determined. This value minimizes the deviations of the volume repair works of the first half of the year from the volume of repair works of the second half the year. Then similarly $x_3$ and $x_9$ are determined. After determination of $x_3,x_6,x_9$ the optimal distribution of the repair works volume is produced for each quarter.

\textbf{Method 3.} Let us suppose $A_1,A_2,\dots,A_{12}$. Let us start from the first month. If $A_1=A_{avg}$ in the first month, then  we go to the second month. If $A_1 < A_{avg}$, then we subtract the missing hours from the second month. If we have an excess of repair works in the first month, then we transfer them to the second month. As a result of using this algorithm we have the following: the transfer from the $i$-th to the $i+1$-th month is denoted by $x>0$, the lack of transfers is denoted by $x=0$, transfer to the $i$-th month from $i+1$ is denoted by $x<0$.

Let us describe the algorithm by steps.

Step 1. Let us suppose $I=1$.

Step 2. Let us calculate $X=A_i-A_{avg}$. Then let us consider the $i$-th month:
a) if $A_i > A_{avg}$, then the value of repair works equals $A_{avg}$ in the $i$-th month. Then we transfer to the $i+1$-month the rest of the repair works: $A_i=A_i-X_i$. Then the step 4 is executed.
b) if $A_i=A_{avg}$, then nothing is transferred from the $i$-th month.  Then the step 4 is executed.
c) if $A_i < A_{avg}$, then we execute to the step 3.

Step 3. Let us compare $A_{i+1}$ and $(-X_i)$:
a) if $A_{i+1}\geq -X_i$, then $A_i=A_i-X_i$; $A_{i+1}=A_{i+1}+X_i$. Then we execute the step 4.
b) if $A_{i+1}< X_i$, then $X_i=A_{i+1}$; $A_{i}=A_{i}-X_i$, $A_{i+1}=A_{i+1}+X_i$. Then we execute the step 4.

Step 4. Let us suppose $I = I + 1$.

Step 5. If $I\geq 12$, then we execute the step 2, otherwise we execute the step 6.

Step 6. As a result we have $X_i,A_j,i=\overline{1,11}, j=\overline{1,12}$.

Step 7. Algorithm stops.

\textit{The solution of the problem 2.}

Let us suppose, that the monthly transfers $X={X_1,\dots,X_{11}}$ are found using one of the three methods. With the help of simple transformation we get the following sets consisting of the elements of original set.
The first set $z=\{z_1,\dots,z_{11}\}$, where $z_i>0$ if $z_{ij}\in X$. Otherwise $z_i=0$. If $z_i>0$, then we transfer from $i$-th month to $i+1$-th month.
The second set $y=\{y_1,\dots,y_{12}\}$, where $y_i>0$ if $y_i\in X$ (transfer from $j$-th month to $j+1$-th month). Otherwise $y_i=0$ (nothing is transferred from $j$-th month to $j+1$-th month).
The goal of introducing the sets $Y$ and $Z$ is the transformation of the problem to any convenient type suitable for application of the methods of dynamic programming.
Let us suppose $D$ is the volume of repair works that is transferred. The vector of monthly repair works is $\delta_n=\{0,1\}$. This vector consists of the numbers $a_{ij}$, that  we now denote by $a_i$. Let us consider the problem $(D-\sum\limits_{i=1}^Na_i\delta_i)=\max\limits_{\{\delta_i\}}\sum\limits_{i=1}^N (a_i\delta_i-D)$=$\max\limits_{\{\delta_i\}}\sum\limits_{i=1}^N(a_i\delta_i-D/N)=\max\limits_{\{\delta_i\}}\sum\limits_{i=1}q_i(\delta_i)$, 
where $q_i(\delta_i)=a_i\delta_i-D/N$ under the constraints $\sum\limits_{i=1}^N \delta_i a_i \leq D,  \delta_i=\{0,1\}$.
Let us write a recurrent relation of dynamic programming in the form of Bellman function: $f_{\overline{N}}(D)=\max\limits_{\{\delta_i\}}\sum\limits_{i=1}^{\overline{N}}q_i(\delta_i)$, $f_N(D)=\max\{\delta_N a_N + f_{N-1}(D-\delta_N a_N)\}$, $\delta_N=0.1$, $\delta_Na_N \leq D$, $N=1,\dots,\overline{N}$.
The above algorithm was implemented for the matrix
\[ \left( \begin{array}{cccc}
10&20 & 30 & 40\\
5 & 8& 6 & 6\\
21 & 11 & 3 & 2\\
14 & 1 & 5 & 3 \end{array} \right).	
\textnormal{\qquad The solution has a form: } 
\left( \begin{array}{cccc}
0&0&0&0\\
0&-1&0&-1\\
1&-1&-1&0\\
0&1&0&0\\
\end{array} \right).\]

\section{Conclusion}
The problem of optimal preventive equipment repair planning is formalized and investigated in this paper. The optimal annual plan provides a uniform distribution of repair works by the months. An algorithm for finding the optimal plan of repair works is described. A numerical example is calculated.

	\section{Acknowledgements}
	
	The authors express gratitude to M. Badalov for solution of the numerical example.
	
	The work is partly supported by work RFBR No. 18-01-00796.


\begin{thebibliography}{}
	
\bibitem{HMWagner} H. M.Wagner, R. J.Giglio, and R. Glaser, Management Science 10, 1964,316–334.
\bibitem{RDBaker}R. D. Baker, The Journal of the Operational Research Society 42, 1964, 493–503.
	

\bibitem{Zub}Malafeev ~ O. \: A., Zubova ~ A. \: F. Mathematical and computer modeling of socio-economic systems at the level of multi-agent interaction (Introduction to the problems of equilibrium, stability and reliability). SPb.:Publishing SPbSU, 2006. P.~1006.

\bibitem{Sosn}Malafeev ~ O. \: A., Sosnina ~ V. \: V. Management model process of cooperative three-agent interaction. Problems of mechanics and control: nonlinear dynamical systems, 2007. №39, p.~131-144.

\bibitem{Grigor}Grigorieva ~ K. \: V., Malafeev ~ O. \: A. The dynamic process of cooperative interaction in the multi-criteria (multi-agent) postman task. The Bulletin of Civil Engineers, 2011. №1, p.~150-156.

\bibitem{Malaf}Malafeev ~ O. \: A. Managed conflict systems. Saint-Petersburg, 2000. P.~280.

\bibitem{Kolok}Malafeev ~ O. \: A., Kolokoltsov ~ V. \: N.
Understanding game theory. New Jersey, 2010. P.~286.

\bibitem{Zenov}Malafeev ~ O. \: A., Zenovich ~ O. \: S., Sevek ~ V. \: K. Multi-agent interaction in the dynamic problem of managing construction projects. Economic Revival of Russia, 2012. № 1, p.~124-131.

\bibitem{Parsh}Drozdova ~ I. \: V., Malafeev ~ O. \: A., Parshina ~ L. \: G. Efficiency of options for the reconstruction of urban housing. Economic Revival of Russia, 2008. № 3, p.~63-67.

\bibitem{Pach}Malafeev ~ O. \: A., Pachar ~ O. \: V.Dynamic, non-stationary task of investing projects in a competitive environment. Problems of mechanics and control: Nonlinear dynamical systems, 2009. № 41, p.~103-108.

\bibitem{Gord}Gordeev ~ D. \: A., Malafeev ~ O. \: A., Titova ~ N. \: D. Probabilistic and deterministic model of the influence factors on the activities of the organization
to innovate. Economic Revival of Russia, 2011. № 1, p.~ 73-82.

\bibitem{Ivan}Grigorieva ~ K. \: V., Ivanov ~ A. \: S.,  Malafeev ~ O. \: A. Static coalition model of investment of innovative projects. Economic Revival of Russia, 2011. № 4, p.~ 90-98.

\bibitem{Chern}Malafeev ~ O. \: A., Chernych ~ K. \: S. Mathematical modeling of the company's development. Economic Revival of Russia, 2004. №, p.~ 60-66.

\bibitem{Tit}Gordeev ~ D. \: A., Malafeev ~ O. \: A., Titova ~ N. \: D. Stochastic model of decision-making about bringing to market an innovative product. Herald of civil engineers, 2011. № 2, p.~ 161-166.

\bibitem{Kolok2}Kolokoltsov ~ V. \: N., Malafeev ~ O. \: A. Mathematical modeling of systems of competition and cooperation (game theory for all), textbook. Saint-Petersburg, 2012. P.~ 624.

\bibitem{Gric}Gricai ~ K. \: N., Malafeev ~ O. \: A. The problem of competitive management in the model of multi-agent interaction of the auction type. Problems of mechanics and control: nonlinear dynamical systems, 2007. №39, p.~36-45.

\bibitem{Akul}Akulenkova ~ I. \: V., Drozdov ~ G. \: D., Malafeev ~ O. \: A. Problems of reconstruction of housing and communal services of a megacity, monograph. Ministry of Education and Science of the Russian Federation, Federal Agency for Education, St. Petersburg State University of Service and Economics, Saint-Petersburg, 2007. P.~187.

\bibitem{Parf}Parfenov ~ A. \: P., Malafeev ~ O. \: A. Equilibrium and compromise control in network models of multi-agent interaction. Problems of mechanics and control: nonlinear dynamical systems, 2007. №39, p.~154-167.

\bibitem{Troev}Malafeev ~ O. \: A., Troeva ~ M. \: S. A weak solution of Hamilton-Jacobi equation for a differential two-person zero-sum game. In the collection: Preprints of the
Eight International Symposium on Differential Games and Applications, 1998. P.~366-369.

\bibitem{Drozd}Drozdova ~ I. \: V., Malafeev ~ O. \: A., Drozdov ~ G. \: D. Modeling the processes of reconcreting the housing and communal services of a metropolis in a competitive environment, monograph. The Federal Agency for Education, St. Petersburg State University of Architecture and Civil Engineering, Saint-Petersburg, 2008. P.~147.

\bibitem{Boit}Malafeev ~ O. \: A., Boitsov ~ D. \: S.,  Redinskikh ~ N. \: D. Compromise and balance in models of multi-agent management in the corrupt network of society. A young scientist, 2014. № 10 (69), p.~14-17.

\bibitem{Grit}Malafeev ~ O. \: A., Gritsai ~ K. \: N., Redinskikh ~ N. \: D. Competitive management in auction models. Problems of mechanics and control: nonlinear dynamical systems, 2004. № 36, p.~74-82.

\bibitem{Ersh}Ershova ~ T. \: A., Malafeev ~ O. \: A. Conflict management in the model of entering the market. Problems of mechanics and control: nonlinear dynamical systems, 2004. № 36, p.~19-27.

\bibitem{Grig2}Grigorieva ~ K. \: V., Malafeev ~ O. \: A. Methods for solving the dynamic multicriteria mailman problem. Herald of civil engineers, 2011. № 4, p.~156-161.

\bibitem{Troev2}Malafeev ~ O. \: A., Troeva ~ M. \: S. Stability and some numerical methods in conflict-control systems. Yakutsk, 1999. P.~102.

\bibitem{Shkr}Shkrabak ~ V. \: S., Malafeev ~ O. \: A., Skrobach ~ A. \: V., Skrobach ~ V. \: F. Mathematical modeling of processes in agro-industrial production. Saint-Petersburg, 2000. P.~336.

\bibitem{Malaf2}Malafeev ~ O. \: A. On the existence of a critical value of a dynamic game. Bulletin of St. Petersburg University. Series 1. Mathematics. Mechanics. Astronomy. Saint-Petersburg, 1972. № 4, p.~41-46.

\bibitem{Murav}Malafeev ~ O. \: A., Muraviev ~ A. \: I. Mathematical models of conflict situations and their resolution. Volume 1. General theory and all sorts of information, Saint-Petersburg, 2000. P.~283.

\bibitem{Drozd2}Malafeev ~ O. \: A., Drozdov ~ G. \: D. Modeling processes in the system of urban construction management. Volume 1, Saint-Petersburg, 2001. P.~401.

\bibitem{Korol}Malafeev ~ O. \: A., Koroleva ~ O. \: A. The model of corruption in contracting. Proceedings Edited by Smirnov ~ N. \: V., Tamasyan ~ G. \: Sh. In the collection: Management processes and persistence. Proceedings of the XXXIX International Scientific Conference of Post-Graduate Students and Students, 2008. P.~446-449.

\bibitem{Murav2}Malafeev ~ O. \: A., Muraviev ~ A. \: I. Modeling of conflict situations in socio-economic systems. Saint-Petersburg, 1998. P.~317.

\bibitem{Drozd3}Drozdov ~ G. \: D., Malafeev ~ O. \: A. Modeling of multi-agent interaction of insurance processes. Monograph. Ministry of Education and Science of the Russian Federation, St. Petersburg State University of Service and Economics, Saint-Petersburg, 2010.

\bibitem{Zub2}Zubova ~ A. \: F., Malafeev ~ O. \: A. The Lyapunov stability and oscillation in economic models. Saint-Petersburg, 2001. P.~101.

\bibitem{Bure}Bure ~ V. \: M., Malafeev ~ O. \: A. Agreed strategy in the repeated final games of N persons. Bulletin of St. Petersburg University. Series 1. Mathematics. Mechanics. Astronomy. Saint-Petersburg, 1995. № 1, p.~41-46.

\bibitem{Malaf3}Malafeev ~ O. \: A. On the existence of the meaning of the pursuit game. Siberian Journal of Operation Research, 1970. № 5, p.~25-36.

\bibitem{Redin}Malafeev ~ O. \: A., Redinskikh ~ N. \: D., Alferov ~ G. \: V. Electric circuits analogies in economics modeling: corruption networks. Proceedings Edited by: Egorov ~ N. \: V., Ovsyannikov ~ D. \: A., Veremey ~ E. \: I. In proceeding: 2nd International Conference on Emission Electronics (ICEE) Selected papers, 2014. P.~28-32.

\bibitem{Malaf4}Malafeev ~ O. \: A. Conflict-driven processes with many participants. The dissertation author's abstract on competition of a scientific degree of physical and mathematical sciencesh, Leningrad, 1987.

\bibitem{Never}Malafeev ~ O. \: A., Neverova ~ E. \: G., Nemnyugin ~ S. \: A., Alferov ~ G. \: V. Proceedings Edited by: Egorov ~ N. \: V., Ovsyannikov ~ D. \: A., Veremey ~ E. \: I. Multicriteria model of laser radiation control. In proceedings: 2nd International Conference on Emission Electronics (ICEE) Selected papers, 2014. P.~33-37.

\bibitem{Kolok3}Kolokoltsov ~ V. \: N., Malafeev ~ O. \: A. Dynamic competitive systems of multi-agent interaction and their asymptotic behavior (Part II). Herald of civil engineers, 2011. № 1, p.~134-145.

\bibitem{Malaf5}Malafeev ~ O. \: A. Stability of solutions to multicriteria optimization problems and conflict-controlled dynamic processes. Saint-Petersburg, 1990.

\bibitem{Alfer}Malafeev ~ O. \: A., Redinskikh ~ N. \: D., Alferov ~ G. \: V., Smirnova ~ T. \: E.  Corruption in the models of the first price auction. In the collection: Management in marine and aerospace systems (UMAS-2014) 7th Russian multiconference on problems of management: Conference materials. GNC RF OAO “Concern” CNII “Electrical Appliance”, 2014. P.~141-146.

\bibitem{Nem}Malafeev ~ O. \: A., Nemnyugin ~ S. \: A., Alferov ~ G. \: V. Charged perticles beam focusing with uncontrollable changing parameters. Proceedings Edited by: Egorov ~ N. \: V., Ovsyannikov ~ D. \: A., Veremey ~ E. \: I.  In proceedings: 2nd
International Conference on Emission Electronics (ICEE) Selected papers, 2014. P.~25-27.

\bibitem{Smirn}Malafeev ~ O. \: A., Redinskikh ~ N. \: D., Smirnova ~ T. \: E.  Network model of investing projects with corruption. Management processes and sustainability, 2015. V. 2, № 1, p.~659-664.

\bibitem{Pichug}Pichugin ~ Yu. \: A., Malafeev ~ O. \: A. On assessing the risk of bankruptcy of a firm. In the book: Dynamic Systems: Steadiness, Management, Optimization, Theses of reports, 2015. P.~204-206.

\bibitem{Malts}Alferov ~ G. \: V., Malafeev ~ O. \: A., Maltseva ~ A. \: S. Game-Theoretic model of inspection by anti-corruption group. In proceeding: AIP Conference Proceedings, 2015. P.~450009.

\bibitem{Gerch}Malafeev ~ O. \: A., Redinskikh ~ N. \: D., Gerchiu ~ A. \: L. Optimization model for the location of corrupt officials in the network. In the book: the construction and operation of energy-efficient buildings (theory and practice taking into account the corruption factor) (Passivehouse) Kolchedantsev ~ L. \: M., Legalov ~ I. \: N., Badin ~ G. \: M.,  Malafeev ~ O. \: A., Aleksandrov ~ E. \: E., Gerchiu ~ A. \: L., Vasilev ~ U. \: G., collective monograph. Borovichi, 2015. P.~128-140.

\bibitem{Gerch2}Malafeev ~ O. \: A., Redinskikh ~ N. \: D., Gerchiu ~ A. \: L. Project investment model with possible corruption. In the book: the construction and operation of energy-efficient buildings (theory and practice taking into account the corruption factor) (Passivehouse) Kolchedantsev ~ L. \: M., Legalov ~ I. \: N., Badin ~ G. \: M.,  Malafeev ~ O. \: A., Aleksandrov ~ E. \: E., Gerchiu ~ A. \: L., Vasilev ~ U. \: G., collective monograph. Borovichi, 2015. P.~140-146.

\bibitem{Chern2}Malafeev ~ O. \: A., Chernych ~ K. \:S. Mathematical modeling of the company's development. Economic Revival of Russia, 2005. №~2, p.~23.

\bibitem{Murav3}Malafeev ~ O. \: A., Muraviev ~ A. \:I. Mathematical models of conflict situations and their resolution. Saint-Petersburg, 2001. Vol. 2, Mathematical foundations of modeling the processes of competition and conflicts in socio-economic systems, p.~294.

\bibitem{Kefel}Kefeli ~ I. \:F., Malafeev ~ O. \: A.  Mathematical Principles of Global Geopolitics. Saint-Petersburg, 2013. P.~204.

\bibitem{Murav3}Malafeev ~ O. \: A., Redinskikh ~ N. \: D., Parfenov ~ A. \:P., Smirnova ~ T. \: E. Corruption in the models of the first price auction. In the collection: Institutes and the mechanism of innovative development: world experience and Russian practice, 2014. P.~250-253. Collection of scientific articles of the 4th International Scientific and Practical Conference. Managing editor: Gorochov ~ A. \: A. 

	
\end{thebibliography}
\end{document}